\documentclass[a4paper,11pt]{article}

\usepackage{graphicx,amssymb,amsthm}

\columnwidth\textwidth

\def\la{\lambda}
\def\om{\omega}

\def\bq{\begin{equation}}
\def\eq{\end{equation}}
\def\li{\langle}
\def\ri{\rangle}

\def\ss{\subset}
\newtheorem{df}{Definition}
\newtheorem{lemma}{Lemma}
\newtheorem{prop}{Proposition}
\def\kr{{\cal R}}

\begin{document}
\begin{center}
{\LARGE A remark on the symmetries of the Riemann curvature
tensor}
\end{center}
The tetralinear form $R(a,b,c,d)$ defined as $\li R(a,b)c,d\ri$,
where $R$ is a Riemann curvature tensor, satisfies two basic
symmetry properties:
\begin{enumerate}
\item{$R(a,b,c,d)=-R(b,a,c,d)=-R(a,b,d,c)$,}
\item{$R(a,b,c,d)+R(b,c,a,d)+R(c,a,b,d)=0$.}
\end{enumerate}
Out of these one can derive other identities,
e.g.~$R(a,b,c,d)=R(c,d,a,b)$. The purpose of this note is to show
that these properties can be (efficiently) expressed as
$SL(2)$-invariance.

\begin{df}
If $V$ is a vector space then let $\kr(V)\ss (V^*)^{\otimes4}$ be
the subspace of the tetralinear forms satisfying the properties
1.~and 2.
\end{df}

\begin{lemma}
Let $(W,\rho)$ be a finite-dimensional representation of $SL(2)$.
If $w\in W$ is invariant with respect to the subgroup of
upper-triangular matrices $SB(2)\ss SL(2)$, then it is
$SL(2)$-invariant.
\end{lemma}

\begin{proof}The $SB(2)$ invariance means that $w$ is a highest
weight vector of weight 0.
\end{proof}

\begin{prop} Let $U$ be an auxiliary 2-dimensional vector space
with a chosen area form $\om$. Then, if $V$ is a
finite-dimensional vector space, there is a natural isomorphism
\bq\kr(V)\simeq\left(\bigwedge\nolimits^4(U\otimes
V^*)\right)^{SL(U)}.\eq
\end{prop}

\begin{proof}
Let $u_1$, $u_2$ be a basis of $U$ such that $\om(u_1,u_2)=1$ (the
isomorphism will be constructed using this choice, but due to the
$SL(U)$-invariance, it will not depend on it). It gives us an
isomorphism $U\otimes V^*\simeq (u_1\otimes V^*)\oplus (u_2\otimes
V^*)$, i.e. \bq\bigwedge\nolimits^4(U\otimes
V^*)\simeq\bigoplus_{i=0}^4\bigwedge\nolimits^i(u_1\otimes
V^*)\otimes\bigwedge\nolimits^{4-i}(u_2\otimes V^*).\eq If we
restrict the left hand side to the elements invariant w.r.t.~the
diagonal matrices in $SL(2)$ (acting by $u_1\mapsto\mu u_1$,
$u_2\mapsto\mu^{-1} u_2$), we get $\bigwedge\nolimits^2(u_1\otimes
V^*)\otimes\bigwedge\nolimits^2(u_2\otimes V^*)$; this space can
be identified with the tetralinear forms on $V$ satisfying 1.
Restricting further to the elements invariant w.r.t.~$u_1\mapsto
u_1$, $u_2\mapsto u_2+\la u_1$ we get the forms satisfying both
1.~and 2. We have passed to the subspace of $SB(2)$-invariants; by
Lemma 1 this is the same as $SL(2)$-invariants.
\end{proof}

\noindent {\bf Final remarks.} As a simple example, the identity
$R(a,b,c,d)=R(c,d,a,b)$ follows from invariance w.r.t.~the matrix
$u_1\mapsto u_2$, $u_2\mapsto -u_1$. More importantly, there is a
supergeometry behind the ad hoc construction of this note, making
it more natural; it will (hopefully) appear elsewhere.

\begin{flushright}
Pavol \v Severa\\\it
Dept.~of Theor.~Physics\\
Comenius University\\
Bratislava, Slovakia
\end{flushright}

\end{document}